\documentclass{article}
\usepackage{ijcai19}

\usepackage{times}
\usepackage{soul}
\usepackage{url}
\usepackage[hidelinks]{hyperref}
\usepackage[utf8]{inputenc}
\usepackage[small]{caption}
\usepackage{graphicx}
\usepackage{amsmath}
\usepackage{booktabs}
\usepackage{algorithm}
\usepackage{algorithmic}
\usepackage{amsmath}
\usepackage{amssymb}

\newcommand{\Min}{\textrm{minimize}}

\usepackage{theorem}
\newtheorem{lemma}{Lemma}

\newtheorem{theorem}{Theorem}

\newtheorem{definition}{Definition}
\newtheorem{assumption}{Assumption}

\title{Heavy-ball Algorithms Always Escape Saddle Points}

\author{
Tao Sun$^1$\and
Dongsheng Li$^{1}$\footnote{Corresponding Author.}\and
 Zhe Quan$^2$\and Hao Jiang$^1$\and Shengguo Li$^1$ \And Yong Dou$^1$
\affiliations
$^1$ Computer of College, National University of Defense Technology\\
$^2$ Hunan University\\
\emails
  nudtsuntao@163.com,
dsli@nudt.edu.cn,
quanzhe@hnu.edu.cn,\{haojiang, nudtlsg, yongdou\}@nudt.edu.cn
}

\begin{document}

\maketitle

\begin{abstract}
Nonconvex optimization algorithms with random initialization have attracted increasing attention recently. It has been showed that many first-order methods always avoid saddle points with random starting points.
In this paper, we answer a question: can the  nonconvex heavy-ball algorithms with random initialization avoid saddle points? The answer is yes! Direct using the existing proof technique for the  heavy-ball algorithms  is hard due to that each iteration of the heavy-ball algorithm consists of current and last points. It is impossible to formulate the algorithms as iteration like $x^{k+1}=g(x^k)$ under some mapping $g$. To this end, we design a new mapping on a new space. With some transfers, the heavy-ball algorithm  can be interpreted as iterations after this mapping.  Theoretically, we prove that heavy-ball gradient descent    enjoys larger stepsize than the gradient descent to escape saddle points to escape the saddle point. And the  heavy-ball proximal point algorithm is also considered; we also proved that the algorithm can always escape the saddle point.
\end{abstract}

\section{Introduction}\label{intro}
We consider the smooth nonconvex  optimization problem
\begin{align}\label{model}
\Min_{x\in \mathbb{R}^N} f(x)
\end{align}
where $f$ is a differentiable closed function whose gradient is $L$-Lipschitz continuous (may be nonconvex). And more we assume $f$ is $\mathcal{C}^2$ over $\mathbb{R}^N$, that is, $f$ is twice-differentiable and $\nabla^2 f(x)$ is continuous. This paper is devoted to the study of two heavy-ball algorithms for minimizing $f$  with random initializations. Specifically, we focus on whether these two algorithms escape the saddle points.

\subsection{Heavy-ball Algorithms}
 We consider the Heavy-Ball Gradient Descent (HBGD)   with random initializations $x^0,x^1$  presented as
\begin{align}\label{hb}
x^{k+1}=x^k-\gamma \nabla f(x^k)+\beta(x^k-x^{k-1}), k\geq 1,
\end{align}
where $\gamma$ is the stepsize and $\beta$ is the inertial parameter. If $\beta=0$, the algorithm then reduces to the basic gradient descent.
The heavy-ball method  is also named as momentum method and has been widely used in different cases to accelerate the gradient descent iteration. Theoretically, it has been proved that  HBGD enjoys  better convergence factor than both the gradient and Nesterov’s accelerated
gradient method with linear convergence
rates under the condition that the objective function is twice continuously differentiable, strongly convex and
has Lipschitz continuous gradient.
With the convex and smooth assumption on the objective function,  the ergodic $O(1/k)$ rate in terms of the objective value, i.e., $f\left(\frac{\sum_{i=1}^k x^i}{k}\right)-\min f=O\left(\frac{1}{k}\right)$, is proved by \cite{ghadimi2015global}. HBGD was proved to converge linearly in the strongly convex case  by \cite{ghadimi2015global}. But the authors used a somehow restrictive assumption on the parameter $\beta$, which leads to a small range values of choose for  $\beta$
 when the strongly convex constant is tiny.
 Due to the inertial term $\beta(x^k-x^{k-1})$,
 HBGD breaks Fej\'{e}r monotonicity that  gradient descent obeys. Consequently, it is difficult to prove its non-ergodic convergence rates. To this end, \cite{sun2018non} designed a novel Lyapunov analysis and derived better convergence results like non-ergodic  sublinear convergence rate $\left(f(x^k)-\min f=O\left(\frac{1}{k}\right)\right)$ and larger stepsize for linear convergence.  \cite{Combettes2017Quasinonexpansive} developed and studied heavy-ball methods for operator research.
In the nonconvex case, if the objective function is nonconvex but Lipschitz differentiable, \cite{zavriev1993heavy} proved   that the sequence generated by the heavy-ball gradient descent  is  convergent to a critical point, but without specifying related convergence rates. With semi-algebraic assumption on the objective function, \cite{ochs2014ipiano} proved the sequence  generated by   HBGD is convergent to some critical point. Actually,  \cite{ochs2014ipiano} also extended HBGD to a nonsmooth composite optimization problem.  \cite{xu2013block} developed and analyzed the Heavy-ball algorithm for tensor minimization problems. \cite{loizou2017linearly,loizou2017momentum} introduced the stochastic versions of heavy-ball algorithms and its relations with existing methods.

On the other hand, \cite{alvarez2000minimizing} studied a second-order ODE whose  implicit discretization yields
the Heavy-Ball Proximal Point Algorithm (HBPPA) mathematically described as
\begin{equation}\label{ppa}
    x^{k+1}=\textbf{Prox}_{\gamma f}(x^k+\beta(x^k-x^{k-1})),
\end{equation}
which admits the acceleration of  the  proximal point algorithm. In fact, the setting $\beta=0$ will push the scheme back to be the basic  proximal point algorithm.
The convergence   of HBPPA is studied
by \cite{alvarez2000minimizing} with  several assumptions on the parameters $\gamma$ and $\beta$ when $f$ is convex.
Noticing the fact that $\nabla f$ is actually the maximal monotone
operator in the convex case, \cite{alvarez2001inertial} extended HBPPA to solve the inclusion problem.
The inexact HBPPA was proposed and studied by \cite{moudafi2003approximate}. Existing works on HBPPA are all based on the convex assumption on the $f$. In this paper, we will present the convergence of nonconvex  HBPPA including the gradient convergence under smooth assumption and sequence convergence under semi-algebraic assumption.

The performances of  HBGD and  HBPPA with convex settings exclude the interest of  this paper. What we care is that can these two algorithms escape the saddle points  when $f$ is nonconvex? This paper answers this problem.

\subsection{First-order Algorithms Escape Saddle Points}
In the community of nonconvex continuous optimization,  saddle points, which are usually unavoidable and   outnumber the cardinality of local minima, have long been regarded as a major obstacle. In many applications like tensor decomposition \cite{ge2015escaping}, matrix completion \cite{ge2016matrix} and dictionary learning \cite{sun2017complete}, all local optima are also globally optimal. In some networks, under several conditions,  local optimas are close to the global optima \cite{chaudhari2016entropy}.  Although these conditions are barely satisfied by many networks, the conclusion coincides with existing experimental observations.   Thus if  all local minimizers are relatively good, we just need to push the iterations to escape the saddle points.  And how to escape the saddle points has been a hot topic recently for  first-order optimization algorithms.  \cite{pemantle1990nonconvergence} established the  convergence of the Robbins-Monro stochastic approximation to local minimizers by introducing sufficiently large unbiased noise.
For tensor problem, \cite{ge2015escaping}
provided the quantitative rates on the
convergence of noisy gradient descent. For general smooth nonconvex functions, the analysis of noisy gradient descent is established by \cite{jin2017escape}.  In a latter paper, \cite{jin2017accelerated} considered the noisy accelerated gradient descent and established related convergence theory. Another technique to escape the saddle points is using random initializations.
It has been   shown that, even
without any random perturbations, a good initialization can make gradient descent
converge  to the global minimum for various non-convex problems like matrix factorization \cite{keshavan2010matrix} , phase retrieval \cite{cai2016optimal}, dictionary learning \cite{arora2015simple}.
\cite{lee2016gradient,lee2017first} proved that for plenty of first-order nonconvex optimization algorithms, with random initializations, they can converge to local minima with probability 1. In the perspective of mathematics, the core techniques consist of two parts: 1). they reformulate the algorithms as mapping iterations; 2). the  stable manifold theorem \cite{shub2013global} is employed. In this paper, we establish the theoretical guarantees for heavy-ball algorithms escaping the saddle points.
\subsection{Difficulty in Analyzing Heavy-ball Iterations}
The difficulty of study on the performance of heavy-ball with random initializations lies in how to reformulate the algorithm   as iterations under a mapping.
Taking the gradient descent for  example, one can use $g:=\mathbb{I}-\gamma\nabla f$ with $\gamma$ being the stepsize. And then, the gradient descent can be represented as $x^{k+1}=g(x^k)$. Constructing the mapping for gradient descent is obvious and easy.    But for heavy-ball ones, the construction tasks are quite different due to the inertial term $\beta(x^k-x^{k-1})$. Point $x^{k+1}$ is generated by both the current variable $x^k$ and the last one $x^{k-1}$. Thus, it is impossible to find a mapping  $g$ to rewrite the heavy-ball as $x^{k+1}=g(x^k)$.
\subsection{Contribution}
Due to that the analysis of nonconvex HBPPA is still blank, we establish  related theory in this paper. And with the recently developed techniques in variational analysis, the sequence convergence of nonconvex HBPPA is also studied.
But our main interests are still
the analysis of HBGD and HBPPA with random initializations. We prove guarantees for both algorithms to escape saddle points under mild assumptions. Direct using the existing method is unapplicable for heavy-ball algorithms. Thus, we propose a novel construction technique. It can be proved that heavy-ball gradient descent enjoys larger stepsize than gradient descent to escape the saddle points with random starting points, which somehow indicates the advantages of heavy-ball methods.

\section{Preliminaries}
In this part, we introduce several basic definitions and tool lemmas. For a closed function $J$ (not necessary to be convex), the proximal map is defined as
\begin{align}\label{prox}
\textbf{Prox}_{J}(\cdot):=\textrm{arg}\min_{x}\left\{J(x)+\frac{\|x-\cdot\|^2}{2}\right\}.
\end{align}
Under the condition that $J$ is differentiable, if $y=\textbf{Prox}_{J}(z)$, the K.K.T. condition can yields that
$\nabla J(y)+y-z=\textbf{0}.$ Given another point arbitrary $x\in \mathbb{R}^N$, it holds that
\begin{equation}\label{proxr}
   J(y)+\frac{\|y-z\|^2}{2}\leq J(x)+\frac{\|x-z\|^2}{2}.
\end{equation}
For matrices $A$ and $B$, $A\sim B$ denotes that there exist  perturbation matrices $P,Q$ such that $A=PBQ$. We use $\mathbb{I}_{N}$ and $\textbf{0}_{N}$  to denote the $N\times N$ unit  and zero matrix, respectively. $\textrm{det}(A)$ is  determinant of matrix $A$. It is easy to see if $A\sim B$, $\textrm{det}(A)=0\Leftrightarrow \textrm{det}(B)=0$.
\begin{definition}[Strict Saddle]\label{saddle}
We present the following two definitions,
\begin{enumerate}
\item A point $x^*$ is a critical point of $f$ if $\nabla f(x^*)=0$. $x^*$
is isolated if there is a neighborhood $U$ around $x^*$, and $x^*$
is the only critical point in $U$.

\item A point $x^*$ is a strict saddle point of $f$ if $x^*$ is a critical point and  $\lambda_{\min} (\nabla^2 f(x^*)) <0$. Let $\chi^*_f$ denote the set of strict saddle points of function $f$.
\end{enumerate}
\end{definition}

Let $g$ be a mapping; and its domain is defined as $\textrm{dom}(g):=\{x\mid g(x)\neq\emptyset\}$. And $g^k$ denotes the $k$-fold composition of $g$. Our interest is   the sequence generated by
\begin{align}\label{scheme}
x^{k}=g(x^{k-1})=g^k (x_0).
\end{align}

\begin{definition}[Global Stable Set]
The global stable set of the strict saddles is the set of initial conditions where iteration of the mapping $g$ converges to a strict saddle, i.e.,
$
W_g := \{ x_0: \lim_k g^k (x_0) \in \chi^* \}.
$
For an iterative algorithm $\mathcal{M}$, $W_{\mathcal{M}}:= \{ x_0: \lim_k x^k  \in \chi^* \}$ and $\{x^1,\ldots\}$ is generated by $\mathcal{M}$ with starting point $x^0$.
\end{definition}

\begin{definition}[Unstable fixed point]
The differential of the mapping $g$, denoted as $\textrm{D} g (x)$, is a linear operator from $\mathcal{T}(x)  \to \mathcal{T}(g(x))$, where $\mathcal{T}(x)$ is the tangent space of $\mathcal{X}$ at point $x$. Given a curve $\gamma$ in $\mathcal{X}$ with $\gamma(0) =x$ and $\frac{d\gamma}{dt}(0):=v \in \mathcal{T}(x)$,  the linear operator is defined as $\textrm{D}g (x)  v = \frac{d( g \circ \gamma) }{dt} (0) \in \mathcal{T}(g(x))$.
Let
$$\mathcal{A}^* _g= \{x: g(x) =x, \exists \lambda (Dg(x))\textrm{ such that }|\lambda (Dg(x))|>1 \}$$ be the set of fixed points where the differential has at least a single real eigenvalue with magnitude greater than one. These are the unstable fixed points.
\end{definition}

With previous definitions, the authors in \cite{lee2017first} presented the following result. The core technique is employing the stable manifold theorem given in \cite{shub2013global}.
\begin{lemma}\label{core}
Let $g$ be a $\mathcal{C}^1$ mapping  and $\det( D g(x))\neq \mathbf{0}$ for all $x \in \emph{dom}(g)$. Assume the sequence $(x^k)_{k\geq 0}$ is generated by scheme \eqref{scheme}. Then the set of initial points that converge to an unstable fixed point has measure zero,  $\mu(\{x^0: \lim x^k \in \mathcal{A}^*_g \} ) =0$. Further if $\mathcal{X}^* \subset \mathcal{A}^*_g$,  $\mu(W_g) =0$.
\end{lemma}

\section{Main Results}
In this section, we present the theoretical guarantees to escape the saddle points for the heavy-ball algorithms presented before. We need a smooth assumption on the objective function.
\begin{assumption}\label{ass}
The objective function $f$ is twice-differentiable and $\|\nabla^2 f(x)\|_2\leq L$, where $\|\cdot\|_2$ denotes the spectral norm of a matrix. Otherwise the objective function satisfies that $\lim_{x\in \mathbb{R}^N} f(x)>-\infty$.
\end{assumption}
The smooth assumption is used to derive the differential of the designed mapping. And the lower bounded assumption is used to provide a  lower bound for the descent variables.

\begin{lemma}\label{th1}
Let $(x^k)_{k\geq 0}$ be generated by the HBGD \eqref{hb} with $0<\beta<1$. If the stepsize is selected as
$$0<\gamma<\frac{2(1-\beta)}{L},$$
we have $\mu\{(x^0,x^1):\lim_{k}x^k\in\chi^*_{f}\}=0$ and $\mu\{W_{\textit{HBGD}}\}=0$.
\end{lemma}

It is easy to see that if $\beta<\frac{1}{2}$ the upper bound of the stepsize is therefore larger than $\frac{1}{L}$ for HBGD with random initializations. And the upper bound is closed to $\frac{2}{L}$ provided that $\beta$ is sufficiently small.   \cite{lee2016gradient}  proved that the stepsize for the GD with random initializations is required to be smaller than $\frac{1}{L}$. Lemma \ref{th1} means a larger stepsize can be used under the heavy-ball scheme, which somehow demonstrates the advantage of heavy-ball method.

The difficulty of the analysis for HBGD has been presented in Sec. 1: it is impossible to find to a mapping such that $x^{k+1}=g(x^k)$ directly. Noticing the $k$th iteration involves with $x^{k}$ and $x^{k-1}$, we turn to consider combining the these two points together, i.e., $(x^{k}, x^{k-1})$. And then, the task  then turns to find a mapping such that $(x^{k+1}, x^{k})=\hat{g}(x^{k}, x^{k-1})$. The mapping $\hat{g}$ is from $\mathbb{R}^{2N}$ to $\mathbb{R}^{2N}$. After building the relations between the unstable fixed point of $\hat{g}$ and the strict saddle points, the theorem then can be proved. Different from existing results in  [Lemma 1,\cite{sun2018non}], $\beta\neq 0$ here.
This is to make sure that
the differential of the
constructed mapping in the proofs is invertible. And the differential is non-symmetric.

 Lemma \ref{th1} actually means $\mathbf{Prob}(x^0\in W_{\textit{HBGD}})=0$, that is, HBGD always escapes saddle points provided the stepsize is set properly. Lemma \ref{th1} just means that $(x^k)_{k\geq 0}$ will not fall into a saddle point. But whether the sequence converges or not  is out of scope of this lemma. In the nonconvex optimization community, such a problem is answered by posing the semi-algebracity assumption on the objective function \cite{lojasiewicz1993geometrie,kurdyka1998gradients} or the isolated assumption on the critical point \cite{lange2013elementary}. Thus, we can derive the following result.
\begin{theorem}\label{pro1}
Assume that conditions of Lemma \ref{th1} hold and $f$ is coercive\footnote{We say function coercive if $f(x)\rightarrow+\infty$ as $\|x\|\rightarrow+\infty$.}. Suppose the starting points $x^0$ and $x^1$ are totally random and one of the following conditions is satisfied.
\begin{enumerate}
  \item $f$ is semi-algebraic.
  \item  All critical points of $f$ are isolated points.
\end{enumerate}
Then, for HBGD, $\lim_{k} x^k$ exists and that limit $x^*$ is a local minimizer of $f$ with probability 1.
\end{theorem}

Now, we study the HBPPA. First, we provide the convergence in the nonconvex settings.
\begin{lemma}\label{leconver}
Let $(x^k)_{k\geq 0}$ be generated by scheme \eqref{ppa} with $0<\beta<\frac{1}{2}$ and $\gamma>0$.
Then, we have $\lim_{k}\|\nabla f(x^k)\|=0$. Further more, if the function $f$ is semi-algebraic and coercive, $(x^k)_{k\geq 0}$  is
convergent to some critical point $x^*$.
\end{lemma}
To obtain the sequence convergence of HBPPA, besides the semi-algebraicity, the coercivity is also needed, which is used to guarantee the boundedness of the generated sequence.
In Lemma \ref{leconver}, $f$ is set to be differentiable. Actually, the smooth setting can be removed. Alternatively, we need to use the notion about subdifferentials of closed functions \cite{rockafellar2009variational} which are frequently used in the variational  analysis. Specifically, in the smooth case, subdifferential is actually the gradient.   With the same parameters, the sequence of HBPPA is still provably to converge to a critical point of $f$.

We are prepared to present the guarantee for HBPPA to escape the saddle point.
\begin{lemma}\label{th2}
Let $(x^k)_{k\geq 0}$ be generated by scheme \eqref{ppa} with $0<\beta<\frac{1}{2}$,
if the stepsize is selected as
$$0<\gamma<\frac{1}{L},$$
we have $\mu\{(x^0,x^1):\lim_{k}x^k\in\chi^*_{f}\}=0$ and $\mu\{W_{\textit{HBPPA}}\}=0$.
\end{lemma}
For HBPPA, the stepsize requirement coincides with existing result \cite{lee2017first}. Main idea of proof of
Lemma \ref{th2} is similar to Lemma \ref{th1}, but more complicated in details.

\begin{theorem}\label{pro2}
Assume that conditions of Lemma \ref{th2} hold and $f$ is coercive. Suppose the starting points $x^0$ and $x^1$ are totally random and one of the following conditions is satisfied.
\begin{enumerate}
  \item $f$ is semi-algebraic.
  \item  All critical points of $f$ are isolated points.
\end{enumerate}
Then, for HBPPA, $\lim_{k} x^k$ exists and that limit $x^*$ is a local minimizer of $f$ with probability 1.
\end{theorem}
\section{Proofs}
This section collects the proofs for the this paper.
\subsection{Proof of Lemma \ref{th1}}
To guarantee the convergence of HBGD, from [Lemma 1,\cite{sun2018non}], we need $0\leq\beta<1$ and $0<\gamma<\frac{2(1-\beta)}{L}$.
Let $y,z\in \mathbb{R}^N$ and $w:=(y,z)$, denote a map as
\begin{align*}
F(w)=F(y,z):=\left(
                \begin{array}{c}
                  y-\gamma\nabla f(y)+\beta(y-z) \\
                  y\\
                \end{array}
              \right).
\end{align*}
It is easy to check that if $F(w)=\textbf{0}\Longrightarrow y=z$.
We further denote
$w^k:=(x^k,x^{k-1})\in \mathbb{R}^{2N}.
$
It can be easily verified that the Heavy-ball methods can be reformulated as
\begin{align}\label{grad}
w^{k+1}= F(w^k).
\end{align}
To use Lemma \ref{core}, we turn to the proofs of  two facts:
\begin{enumerate}
\item $\det( D F(w))\neq \mathbf{0}$ for all $w \in \mathbb{R}^{2N}$.

\item $\chi^*_f\times \chi^*_f\bigcap\{(y,z)\in\mathbb{R}^{2N}\mid y=z\}\subseteq \mathcal{A}_{F}^*$.
\end{enumerate}

Proof of fact 1: Direct calculations give us
\begin{equation*}
    D F(w)=\left(
             \begin{array}{cc}
               (1+\beta)\mathbb{I}_{N}-\gamma\nabla^2 f(y) & -\beta \mathbb{I}_{N}\\
               \mathbb{I}_{N} & \mathbf{0}_{N} \\
             \end{array}
           \right).
\end{equation*}
We use the short-hand notation
$A(y):=(1+\beta)\mathbb{I}_{N}-\gamma\nabla^2 f(y).$
We can obtain the following relation
$$D F(w)=\left(
             \begin{array}{cc}
              A(y) & -\beta \mathbb{I}_{N} \\
               \mathbb{I}_{N}  &   \textbf{0}_{N} \\
             \end{array}
           \right)\sim\left(
                        \begin{array}{cc}
                          \textbf{0}_{N}  &  -\beta \mathbb{I}_{N} \\
                         \mathbb{I}_{N}  &  \textbf{0}_{N} \\
                        \end{array}
                      \right).
           $$
Therefore, we can derive $\textrm{det}\left(D F(w)\right)\neq 0$ for any $w\in \mathbb{R}^{2N}$.

Proof of fact 2: For any $x^*$ being a   strict saddle, with the symmetry of $\nabla^2 f(x^*)$, then
$$ \lambda_{\max}(A(x^*))>1+\beta.$$
For $\lambda\in \mathbb{R}$ and $\lambda\neq 0$~~\footnote{We have proved that $D F(w)$ is nonsingular.}, denoting $w^*=(x^*,x^*)$, we consider
$$\lambda \mathbb{I}_{2N}-D F(w^*)=\left(
             \begin{array}{cc}
               \lambda \mathbb{I}_{N} -A(x^*) & \beta \mathbb{I}_{N} \\
               -\mathbb{I}_{N}  & \lambda \mathbb{I}_{N} \\
             \end{array}
           \right).$$
After simplifications, we can get
\begin{align*}
&\left(
             \begin{array}{cc}
               \lambda \mathbb{I}_{N}-A(x^*) & \beta \mathbb{I}_{N} \\
               -\mathbb{I}_{N}  & \lambda \mathbb{I}_{N} \\
             \end{array}
           \right)\\
&\quad\quad\sim\left(
                        \begin{array}{cc}
                          \lambda \mathbb{I}_{N}+\frac{\beta}{\lambda} \mathbb{I}_{N}-A(x^*) &  \\
                           &  \lambda \mathbb{I}_{N}\\
                        \end{array}
                      \right).
\end{align*}
That indicates
$\textrm{det}(\lambda \mathbb{I}_{2N}-D F(w^*))=0\Longleftrightarrow \textrm{det}(\lambda \mathbb{I}_{N}+\frac{\beta}{\lambda} \mathbb{I}_{N}-A(x^*))=0.$
Denote function $s(t):=t+\frac{\beta}{t},~t>0$. All the eigenvalues of $D F(w^*)$ are the roots of the equation above. We exploit the two roots of
$\lambda+\frac{\beta}{\lambda}=\lambda_{\max}(A(x^*)). $
It is easy to check that
$\lambda_{\max}(A(x^*))^2-4\beta>(1-\beta)^2\geq 0,$
which means this equation
has two real roots
denoted as $0<\underline{\lambda}<\overline{\lambda}$.
Noting $\lambda_{\max}(A(x^*))>1+\beta$, that is
$s(\overline{\lambda})>s(1).$
We can see that $s(t)$ is increasing on $[\sqrt{\beta},+\infty)$.
With the fact $0\leq \beta<1$, $s(t)$ is increasing on $[1,+\infty)$. Thus, there exists some  real number  $\lambda_0>1$  such that $\overline{\lambda}>\lambda_0$. Thus, we get
$\lambda_{\max}(D F(w^*))\geq\overline{\lambda}\geq \lambda_0>1,$
i.e., $w^*\in \mathcal{A}_{F}^*$. That is also
$$\{(x^0,x^1):\lim_{k}x^k\in\chi^*_{f}\}\subseteq\{w^0:\lim_{k}w^k\in \mathcal{A}_{F}^*\}.$$
Using Lemma \ref{core}, $\mu\{w^0:\lim_{k}w^k\in \mathcal{A}_{F}^*\}=0$ and then
$$\mu\{(x^0,x^1):\lim_{k}x^k\in\chi^*_{f}\}\leq\mu\{w^0:\lim_{k}w^k\in \mathcal{A}_{F}^*\}=0.$$
On the other hand, due to the nonnegativity of the measurement,
$$\mu\{(x^0,x^1):\lim_{k}x^k\in\chi^*_{f}\}=0.$$

\subsection{Proof of Theorem \ref{pro1}}
The first case has been proved in  \cite{ochs2014ipiano}, in which  $(x^k)_{k\geq 0}$ is convergent to some critical point $x^*$. For the second case,   [Proposition 12.4.1, \cite{lange2013elementary}] and equation \eqref{leconver-t3} means that the stationary set of $(x^k)_{k\geq 0}$ is finite. Otherwise,  the stationary set is connected which contradicts with the fact all critical points are isolated. Once using [Proposition 12.4.1, \cite{lange2013elementary}], the stationary set has only one element $x^*$; and then, $(x^k)_{k\geq 0}\rightarrow x^*$.

Therefore, in both cases, it can be proved that $(x^k)_{k\geq 0}\rightarrow x^*\in \mathbb{R}^N$. With Lemma \ref{th1}, $x^*$ is a local minimizer with probability 1 if the starting points $x^0,x^1$ are random. The theorem is then proved.
\subsection{Proof of Lemma \ref{leconver}}
This subsection contains two parts. Part 1 is to prove the  gradient convergence, while Part 2 contains the sequence convergence proof.

\textbf{Part 1.}
In \eqref{proxr}, by substituting $y\leftarrow x^{k+1}$, $z\leftarrow x^k+\beta(x^k-x^{k-1})$, $x\leftarrow x^{k}$, $J\leftarrow \gamma f$,
\begin{align*}
   &\gamma f(x^{k+1})+\frac{\|x^k+\beta(x^k-x^{k-1})-x^{k+1}\|^2}{2}\\
   &\quad\quad\leq \gamma f(x^k)+\frac{\|\beta(x^k-x^{k-1})\|^2}{2}.
\end{align*}
After simplifications, we then get
\begin{align}\label{leconver-t1}
&f(x^{k+1})+\frac{\|x^{k+1}-x^k\|^2}{2\gamma}\nonumber\\
&\leq f(x^k)+\frac{\beta}{\gamma}\langle x^{k+1}-x^k,x^k-x^{k-1}\rangle\nonumber\\
&\leq f(x^k)+\frac{\beta}{2\gamma}\| x^{k+1}-x^k\|^2+\frac{\beta}{2\gamma}\|x^k-x^{k-1}\|^2.
\end{align}
Denote a novel function as
$$H(w)=H(x,y):=f(x)+\frac{\beta}{2\gamma}\|x-y\|^2$$
and sequence $w^k:=(x^k,x^{k-1})$.
Thus, inequality \eqref{leconver-t1} offers the descent as
\begin{align}\label{leconver-t2}
H(w^{k+1})+\frac{1-2\beta}{2\gamma}\|x^{k+1}-x^k\|^2\leq H(w^{k}), ~k\geq 1.
\end{align}
With the fact $\inf_{x\in \mathbb{R}^N} f(x)>-\infty$, $\inf_{k} H(w^{k})>-\infty$. The setting $0<\beta<\frac{1}{2}$ gives
\begin{align}\label{leconver-t3}
\lim_k\|x^{k+1}-x^k\|=0.
\end{align}
On the other hand, from the property of the proximal map,
\begin{equation}\label{leconver-t+}
    \gamma\nabla f(x^k)+x^{k+1}-x^k-\beta(x^k-x^{k-1})=\textbf{0},
\end{equation}
which can derive
\begin{align}\label{leconver-t4}
\small
   &\lim_{k} \|\nabla f(x^k)\|\nonumber\\
   &\quad\leq\frac{\lim_{k}\|x^{k+1}-x^k\|+\beta\lim_{k}\|x^k-x^{k-1}\|}{\gamma}=0
\end{align}

\textbf{Part 2.} We prove that $H(x,y)$ is coercive provided $f$ is coercive. If this claim fails to hold, there exists some $C>0$ such that
$$f(x)\leq H(x,y)\leq C~~\textrm{as}~~\|(x,y)\|\rightarrow+\infty.$$
With the coercivity of $f$, there exists some $C_1>0$ such that $\|x\|\leq C_1$. On the other hand, $\|x-y\|$ is bounded due to the lower boundedness of $f$; thus, $\|y\|$ is bounded, which contradicts the fact $\|(x,y)\|\rightarrow+\infty$.

From \eqref{leconver-t2}, we see that $(w^k)_{k\geq 1}$ is bounded; and then, $(x^k)_{k\geq 0}$ is bounded. The descent property of $(H(w^k))_{k\geq 1}$ and the continuity of function $H$  directly give that
\begin{equation}\label{happ}
    \lim_{k}H(w^k)=H(w^*).
\end{equation}
That means  the sequence has a stationary point $x^*$; \eqref{leconver-t4} means $x^*$ admits the critical point of $f$. It is easy to see $w^*:=(x^*,x^*)$ is the critical point of $H$. Without loss of generality, we assume $w^k\neq w^*$ as $k\geq 1$.

Noting polynomial functions are semi-algebraic and the fact that sum of semi-algebraic functions is still semi-algebraic, $H(w)$ is then semi-algebraic. Denote the set of all stationary points of $(w^k)_{k\geq 0}$ as $\Omega$.
From [Lemma 3.6, \cite{bolte2014proximal}],
 there
 exist $\eta,\varepsilon\in (0, +\infty]$ and a continuous concave function $\varphi: [0, \eta)\rightarrow \mathbb{R}^+$ such that
\begin{enumerate}
  \item $\varphi(0)=0$; $\varphi$ is $\mathcal{C}^1$ on $(0, \eta)$;  for all $s\in(0, \eta)$, $\varphi^{'}(s)>0$.
  \item for  $w\in \{w|H(w^*)<H(w)<H(w^*)+\eta\}\bigcap\{w|\textrm{dist}(w,\Omega)<\varepsilon\}$, it holds that
\begin{equation}\label{semi-algebraic}
  \varphi^{'}(H(w)-H(w^*))\cdot \|\nabla H(w)\|\geq 1.
\end{equation}
\end{enumerate}
For the $\eta,\varepsilon$ given above, as $k$ is large enough,  $w^k$ will fall into the set $ \{w|H(w^*)<H(w)<H(w^*)+\eta\}\bigcap\{w|\textrm{dist}(w,\Omega)<\varepsilon\}$.
The concavity of $\varphi$ gives
\begin{align}\label{leconver-t5}
\small
&\varphi(H(w^{k})-H(w^*))-\varphi(H(w^{k+1})-H(w^*))\nonumber\\
&\quad\geq \varphi'(H(w^{k})-H(w^*))\cdot(H(w^{k})-H(w^{k+1}))\nonumber\\
&\quad\overset{a)}{\geq} \nu\cdot\varphi'(H(w^{k})-H(w^*))\cdot\|x^{k+1}-x^{k}\|^2\nonumber\\
&\quad\overset{b)}{\geq} \frac{\nu\|x^{k+1}-x^{k}\|^2}{\|\nabla H(w^k)\|},
\end{align}
where $\nu=\frac{1-2\beta}{2\gamma}$, and  $a)$ uses \ref{leconver-t2}, and $b)$ comes from \eqref{semi-algebraic}. The gradient of $H$ satisfies
\begin{align}\label{leconver-t6}
    &\|\nabla H(w^k)\|\leq \|\nabla f(x^k)\| +\frac{\beta}{\gamma}\|x^k-x^{k-1}\|\nonumber\\
    &\quad\quad\overset{c)}{\leq}\frac{1}{\gamma}(\|x^{k+1}-x^{k}\|+\|x^k-x^{k-1}\|),
\end{align}
where $c)$ depends on \eqref{leconver-t+} and,  we  used the fact $0<\beta<\frac{1}{2}$. Combining  \eqref{leconver-t5} and \eqref{leconver-t6},
\begin{align} \label{leconver-t7}
\small
    &2\|x^{k+1}-x^k\|\nonumber\\
    &\leq2\times \frac{1}{2}\sqrt{\|x^{k+1}-x^{k}\|+\|x^k-x^{k-1}\|}\nonumber\\
    &\times2\sqrt{\frac{1}{\nu\gamma}}\sqrt{\varphi(H(w^{k})-H(w^*))-\varphi(H(w^{k+1})-H(w^*))}\nonumber\\
    &\overset{d)}{\leq}\frac{4}{\nu\gamma}\big[\varphi(H(w^{k})-H(w^*))-\varphi(H(w^{k+1})-H(w^*))\big]\nonumber\\
    &+\frac{\|x^{k+1}-x^{k}\|+\|x^k-x^{k-1}\|}{4},
\end{align}
where $d)$ employs the  Schwarz  inequality $2(xy)^{\frac{1}{2}}\leq x+y$ with $x=2\sqrt{\frac{\zeta}{\nu}}\sqrt{\varphi(H(w^{k})-H(w^*))-\varphi(H(w^{k+1})-H(w^*))}$, and $y= \frac{1}{2}\sqrt{\|x^{k+1}-x^{k}\|+\|x^k-x^{k-1}\|}$. Summing both sides of \eqref{leconver-t7} from sufficiently large $k$ to $K$, we have
\begin{align}\label{leconver-t8}
    &\frac{3}{2}\sum_{i=k}^{K-1}\|x^{i+1}-x^i\|+\frac{7}{4}\|x^{K+1}-x^K\|\leq \frac{1}{4}\|x^{k}-x^{k-1}\|\nonumber\\
    &+\frac{4}{\nu\gamma}\big[\varphi(H(w^{k})-H(w^*))-\varphi(H(w^{K+1})-H(w^*))\big].
\end{align}
Letting $K\rightarrow+\infty$, from \eqref{leconver-t3}, \eqref{happ}, $\varphi(0)=0$ and the continuity of $\varphi$, \eqref{leconver-t8} then yields
$
    \sum_{i=k}^{+\infty}\|x^{i+1}-x^i\|\leq \frac{1}{6}\|x^{k}-x^{k-1}\|+\frac{4}{\nu\gamma}\varphi(H(w^{k})-H(w^*))<+\infty.
$
We are then led to
 $\sum_{i=k}^{+\infty}\|x^{i+1}-x^i\|<+\infty\Rightarrow\sum_{i=0}^{+\infty}\|x^{i+1}-x^i\|<+\infty.$
  That means $(x^k)_{k\geq 0}$ is convergent to some point. With the fact, $x^*$ is a stationary point, $(x^k)_{k\geq 0}\rightarrow x^*$.
\subsection{Proof of Lemma \ref{th2}}
In this proof, we  denote the mapping as
\begin{align*}
F(w)=F(y,z):=\left(
                \begin{array}{c}
                 \textbf{Prox}_{\gamma f}(y+\beta(y-z)) \\
                  y\\
                \end{array}
              \right).
\end{align*}
To calculate the differential of $F$, we denote $$g(y,z):=\textbf{Prox}_{\gamma f}(y+\beta(y-z)).$$ It is easy to see that
\begin{equation}\label{th2-t0}
    D F(w)=\left(
             \begin{array}{cc}
              \frac{\partial g(y,z)}{\partial y} & \frac{\partial g(y,z)}{\partial z}\\
               \mathbb{I}_N  & \mathbf{0}_N \\
             \end{array}
           \right).
\end{equation}
The definition of the proximal map then gives
\begin{equation}\label{th2-t+0}
   \textbf{0}=\gamma\nabla f[g(y,z)]+[g(y,z)-(y+\beta(y-z))].
\end{equation}
By implicit differentiation on variable $y$, we can get
\begin{equation*}
   \textbf{0}=\gamma\nabla^2 f[g(y,z)]\frac{\partial g(y,z)}{\partial y}+\frac{\partial g(y,z)}{\partial y}-(1+\beta)\mathbb{I}_N.
\end{equation*}
Due to that $0<\gamma<\frac{1}{L}$, $\gamma\nabla^2 f[g(y,z)]+\mathbb{I}_N$ is invertible.
We use a shorthand notation
$$A(y,z):=(\gamma\nabla^2 f[g(y,z)]+\mathbb{I}_N)^{-1}.$$
For $z$, by the same procedure, then we can derive
\begin{equation}\label{th2-t1}
   \frac{\partial g(y,z)}{\partial y}=(1+\beta)A(y,z),~\frac{\partial g(y,z)}{\partial z}= -\beta A(y,z).
\end{equation}
With \eqref{th2-t0} and \eqref{th2-t1}, we are then led to
\begin{align}
    D F(w)&= \left(
             \begin{array}{cc}
              (1+\beta)A(y,z) & -\beta A(y,z)\\
               \mathbb{I}_N  & \mathbf{0}_N \\
             \end{array}
           \right)\nonumber\\
           &\sim \left(
             \begin{array}{cc}
              \mathbf{0}_N & -\beta A(y,z)\\
               \mathbb{I}_N  & \mathbf{0}_N \\
             \end{array}
           \right).
\end{align}
That  means for any $w\in \mathbb{R}^{2N}$, $\textrm{det}(D F(w))\neq 0$. On the other hand, if $x^*$ is   a strict saddle, $g(x^*,x^*)=x^*$,
$ \lambda_{\max}(A(x^*,x^*))>1.$
For $\lambda\in \mathbb{R}$, denoting $w^*=(x^*,x^*)$, we consider
\begin{align*}
&\lambda \mathbb{I}_{2N}-D F(w^*)\\
&\quad=\left(
             \begin{array}{cc}
               \lambda \mathbb{I}_{N} -(1+\beta)A(x^*,x^*) & \beta A(x^*,x^*) \\
               -\mathbb{I}_{N}  & \lambda \mathbb{I}_{N} \\
             \end{array}
           \right).
\end{align*}
After simplifications, we can derive
\begin{align*}
&\lambda \mathbb{I}_{2N}-D F(w^*\\
&\quad\sim\left(
                        \begin{array}{cc}
                          \lambda \mathbb{I}_{N}+[\frac{\beta}{\lambda}-(1+\beta)]A(x^*,x^*) &  \\
                           &  \lambda \mathbb{I}_{N}\\
                        \end{array}
                      \right).
\end{align*}
With direct computations, we turn to
$$\textrm{det}(\lambda^2 \mathbb{I}_{N}+(\beta-(1+\beta)\lambda)A(x^*,x^*))=0.$$
We consider the following equation
$$H(\lambda):=\lambda^2+(\beta-(1+\beta)\lambda)\lambda_{\max}(A(x^*,x^*))=0. $$
Direct calculations give
\begin{align*}
&(1+\beta)^2\lambda^2_{\max}(A(x^*,x^*))-4\beta\lambda_{\max}(A(x^*,x^*))\\
&\quad=4\beta\lambda_{\max}(A(x^*,x^*))(\lambda_{\max}(A(x^*,x^*))-1)\\
&\quad+(1-\beta)^2\lambda^2_{\max}(A(x^*,x^*))\geq 0.
\end{align*}
Thus, the equation enjoys two real roots denoted by $\underline{\lambda}<\overline{\lambda}$.
It is easy to see that $\underline{\lambda}>0$. Noticing that $H(1)<0$ and $\lim_{\lambda\rightarrow+\infty}H(\lambda)=+\infty$, we have $\overline{\lambda}>1$. That means
$\lambda_{\max}(D F(w^*))\geq\overline{\lambda}>1,$
 thus, $w^*\in \mathcal{A}_{F}^*$. Consequently, the proof is proved by  Lemma \ref{core}.
\subsection{Proof of Theorem \ref{pro2}}
This proof is similar to the one of Theorem \ref{pro1} and will not be reproduced.

\section{Conclusion}
In this paper, we proved that HBGD and HBPPA always escape the saddle points with random initializations. This paper also established the convergence of nonconvex HBPPA. The core part in the proofs is bundling current and the last point  as one point in an enlarged space. The heavy algorithms then can be represented as an iteration after a mapping. An interesting finding is that the HBGD  can enjoy larger stepsize than the gradient descent to escape saddle points.

\section*{Acknowledgments}
This work is sponsored in
part by National Key R\&D Program of China (2018YFB0204300), and
Major State Research Development Program (2016YFB0201305), and National Natural Science Foundation of Hunan Province
in China (2018JJ3616), and National Natural Science Foundation for the Youth of China (61602166), and Natural Science Foundation of Hunan (806175290082), and Natural Science Foundation of NUDT (ZK18-03-01), and
National Natural Science Foundation of China (11401580).

\end{document}